\documentclass{amsart}
\usepackage{amssymb, amsmath, amsfonts, amscd}

\input xypic

\theoremstyle{plain}
\newtheorem{theorem}{Theorem}[section]
\newtheorem{corollary}[theorem]{Corollary}
\newtheorem{lemma}[theorem]{Lemma}
\newtheorem{proposition}[theorem]{Proposition}

\theoremstyle{definition}
\newtheorem{definition}[theorem]{Definition}

\theoremstyle{remark}

\numberwithin{equation}{theorem}

\renewcommand{\L}{\mathcal{L}}

\renewcommand{\Pr}{\mathcal{P} }

\newcommand{\SL}{\operatorname{SL}} 
 
\newcommand{\Pic}{\operatorname{Pic} }

\newcommand{\End}{\operatorname{End} }

\renewcommand{\H}{\operatorname{H} }

\newcommand{\K}{\operatorname{K} }

\newcommand{\U}{\operatorname{U}} 
\renewcommand{\K}{\operatorname{K} }
 
\renewcommand{\lg}{\mathfrak{g}}
\newcommand{\lp}{\mathfrak{p}}

\newcommand{\ap}{\alpha_P}

\renewcommand{\sl}{\mathfrak{sl}}

\newcommand{\lbp}{\mathfrak{b}_+}

\renewcommand{\ln}{\mathfrak{n}}
\newcommand{\lnm}{\mathfrak{n}_-}
\newcommand{\lnp}{\mathfrak{n}_+}

\newcommand{\lh}{\mathfrak{h}}

\newcommand{\onep}{\mathbf{1}_{\lp}}
\newcommand{\oneg}{\mathbf{1}_{\lg}}

\newcommand{\sym}{\operatorname{Sym} }

\newcommand{\Vl}{V_\lambda}

\newcommand{\ul}{\underline{l} }

\newcommand{\un}{\underline{n} }

\begin{document}

\title{Scalar generalized Verma modules}

\author{Helge Maakestad}
\address{Institut Fourier, Grenoble}

\email{\text{h\_maakestad@hotmail.com}}

\thanks{Supported by a research scholarship from NAV, www.nav.no}


\date{March 2010}

\begin{abstract} In this paper we study the Verma module $M(\mu)$
  associated to a linear form $\mu \in \lh^*$ where $\sl(E)=\lnm\oplus
  \lh \oplus \lnp$ is a triangular decomposition of $\sl(E)$. The
  $\SL(E)$-module $M(\mu)$ has a canonical simple quotient $L(\mu)$
  with a canonical generator $v$. We study the left annihilator ideal
  $ann(v)$ in $\U(\sl(E))$. We also study scalar generalized
  Verma module $M(\rho)$ associated to a character $\rho$ of $\lp$ where $\lp$
  is a parabolic subalgebra of $\sl(E)$. We prove $M(\rho)$ has a
  canonical simple quotient $L(\rho)$. This simple quotient is in some
  cases an infinite dimensional $\SL(E)$-module. We get a class of
  mutually non-isomorphic irreducible $\SL(E)$-modules containing the
  class of all finite dimensional irreducible $\SL(E)$-modules. As a
  result we give an algebraic proof of a classical result of Smoke on
  the structure of the jet bundle as $P$-module on any flag-scheme
  $\SL(E)/P$ where $P$ is any parabolic subgroup.
\end{abstract}

\maketitle

\tableofcontents

\section{Introduction}

We study the scalar generalized
Verma module $M(\rho)$ associated to a character $\rho$ of $\lp$
where $\lp$ is a parabolic subalgebra of $\sl(E)$. We prove $M(\rho)$ has a
canonical simple quotient $L(\rho)$. This simple quotient is in some
cases an infinite dimensional $\SL(E)$-module. We get a construction
of a class of mutually non-isomorphic irreducible $\SL(E)$-modules
containing the class of finite dimensional irreducible
$\SL(E)$-modules. This class contains many infinite dimensional $\SL(E)$-modules.

Let $\Vl$ ba an aritrary finite dimensional irreducible
$\SL(E)$-module.
In this note we give a construction of $\Vl$ using the enveloping
algebra $\U(\sl(E))$ and the Verma modules
$M(\mu)$.  We study the annihilator ideal $ann(v)$ in
$\U(\sl(E))$ where $v$ is the  highest weight vector $v$ in $\Vl$
We prove the class of simple modules $L(\rho)$ may be constructed
using classical Verma modules $M(\mu)$ as done by Dixmier in
\cite{dixmier}. 
As a result we give an algebraic proof of a classical result of Smoke
on the structure of the jet bundle as $P$-module on $\SL(E)/P$ (see
Corollary \ref{jets}).

\section{Scalar generalized Verma modules}

In this section we construct any scalar generalized Verma module
$M(\rho)$. Here $\rho$ is a character of $\lp$ where $\lp$ is any
parabolic subalgebra of $\sl(E)$. We prove $M(\rho)$ has a canonical
simple quotient $L(\rho)$. When $L(\rho)$ is finite dimensional we get
a construction of all finite dimensional irreducible $\SL(E)$-modules.

Let $G=\SL(E)$ where $E$ is an $n$-dimensional vector space over an
algebraically closed 
field $K$ of characteristic zero and let $\lg=\sl(E)$. Let $\lh $ be
the abelian subalgebra of $\sl(E)$ of diagonal matrices. It follows
$(\sl(E),\lh)$ is a split semi simple Lie algebra and determines a
root system $R=R(\sl(E),\lh)$. Let $B$  be a basis for $R$. This
determines the positive roots $R_+$ and the negative roots $R_-$. This
determines a triangular decomposition  $\sl(E)=\lnm\oplus \lh \oplus
\lnm$  of $\sl(E)$.

Let $1\leq n_1<n_2< \cdots < n_k \leq n-1$ be integers where
$dim(E)=n$ and let $l_1,..,l_k \in \mathbf{Z}$. Let $d_1=n_1$,
$d_i=n_i-n_{i-1}$ for $i=2,..,n-1$ and $d_n=n-n_k$.
Let $\ul=(l_1,..,l_k)$ and $\un=(n_1,..,n_k)$.
Let $e_1,..,e_n$ be a basis for $E$ as $K$-vector space  and let
$E_i=K\{e_1,..,e_{n_i}\}$.
Let $P$ in $\SL(E)$ be the subgroup fixing the flag 
\[ E_\bullet: 0\neq E_1\subseteq \cdots \subseteq E_k \subseteq E \]
in $E$. Let $\lp=Lie(P)$. It follows $\lp$ consists of matrices on the
form
\[
x=\begin{pmatrix} A_1 & *      &    \cdots    &    * \\
                  0   & A_2    &    \cdots    &    * \\
                  0   & \vdots &    \cdots    &    * \\ 
                  0   &  0     &    \cdots    & A_{k+1}
\end{pmatrix}
\]
where $A_i$ is a $d_i\times d_i$-matrix and $tr(x)=0$. Let
$L_\rho=Kw$ be a rank one $K$-vector space on the element $w$.
Define the following
character $\rho_{\ul}=\rho$
\[\rho:\lp\rightarrow \End(L_\rho) \]
by
\[ \rho(x)=\eta(x)w=\sum_{i=1}^k l_i(tr(A_1)+\cdots +tr(A_i))w.\]

\begin{lemma} The pair $(L_\rho,\rho)$ is a rank one $\lp$-module.
All rank one $\lp$-modules arise in this way.
\end{lemma}
\begin{proof} The proof is an exercise.
\end{proof}

Let $\U_l(\sl(E))\subseteq \U(\sl(E))$ be the canonical filtration of $\U(\sl(E))$.
\begin{definition} Let $M(\rho)=\U(\sl(E))\otimes_{\U(\lp)}L_\rho$ be
  the \emph{generalized Verma module} associated to $\rho$. Let 
$M_l(\rho)=\U_l(\sl(E))\otimes_{\U(\lp)}L_\rho$
be the \emph{canonical filtration} of $M(\rho)$.
\end{definition}

Since $\U_l(\sl(E))\otimes_{\U(\lp)}L_\rho$ is a $P$-submodule of
$M(\rho)$ where $\lp=Lie(P)$ it follows $\{M_l(\rho)\}_{l\geq 0}$ is a
filtration
of $M(\rho)$ by $P$-modules.

Since $dim_K(L_\rho)=1$ we refer to $M(\rho)$ as a \emph{scalar generalized Verma module}. By definition this construction give all scalar
generalized Verma modules for $\SL(E)$.

Define the following sub Lie algebra of $\sl(E)$: $\ln$ is the
subalgebra of matrices $x$ on the form
\[
x=\begin{pmatrix} A_1 & 0      &    \cdots    &    0 \\
                  *   & A_2    &    \cdots    &    0 \\
                  *   & \vdots &    \cdots    &    0 \\
                  *   &  *     &    \cdots    & A_{k+1}

\end{pmatrix}
\]
where $A_i$ is a $d_i\times d_i$-matrix with zero entries. It follows
there is an isomorphism $\ln\oplus \lp \cong \sl(E)$ as vector spaces.
Let $\sl(E)=\lnm\oplus \lh \oplus \lnp$ be the standard triangular
decomposition of $\sl(E)$ as defined in the previous section. It
follows there is an inclusion $\ln \subseteq \lnm$ of Lie
algebras. Let $R=R(\sl(E),\lh)$ be the roots of $\sl(E)$ with respect
to $\lh$. Let $B'=\{\alpha_1,..,\alpha_m\}$ be a subset of $R$ such
that the set 
\[ X_{-\alpha_1},..,X_{-\alpha_m} \]
is a basis for $\ln$.

Let $P=(p_1,..,p_m)$ and let $X^P=X_{-\alpha_1}^{p_1}\cdots X_{-\alpha_m}^{p_m}$.
Let $\ap=-p_1\alpha_1-\cdots -p_m\alpha_m$.

\begin{lemma} \label{basis} The following holds: The set
\[ \{ X^P\otimes w : p_1,..,p_m\geq 0\} \]
is a basis for $M(\rho)$ as $K$-vector space.
The natural map
\[ \phi:\U(\ln)\rightarrow M(\rho) \]
defined by
\[ \phi(X^P)=X^P\otimes w \]
is an isomorphism of left $\ln$-modules.
\end{lemma}
\begin{proof} Since $\ln\oplus \lp=\sl(E)$ there is by definition an isomorphism
\[ \U(\sl(E))\cong K\{X^P: p_1,..,p_m\geq 0\}\U(\lp) \]
of free right $\U(\lp)$-modules. We get
\[ M(\rho)\cong \{X^P:p_1,..,p_m\geq
0\}\U(\lp)\otimes_{\U(\lp)}L_\rho\cong \]
\[ K\{X^P\otimes w : p_1,..,p_m\geq 0\}.\]
The first claim is proved.
One checks the map $\phi$ is a map of left $\ln$-modules and the Lemma
is proved.
\end{proof}

Let $\omega_i=L_1+\cdots +L_i\in \lh$ for $i=1,..,n-1$ be the fundamental weights for
$\SL(E)$ and let $\lambda=\sum_{i=1}^k l_i\omega_{n_i}$. It follows
for all $x\in \lh \subseteq \lp$ that $\rho(x)=\lambda(x)$.

Let $\rho_{\U}:\U(\lp)\rightarrow \End(L_\rho)$ be the
\emph{associated morphism} of $\rho$. Let $L$ be the following left
ideal of $\U(\lp)$:
\[ L=\U(\lp)\{x-\rho(x)\onep : x\in \lp\}.\]

\begin{proposition} Let $N=ker(\rho_{\U})$. There is an equality of
  ideals in $\U(\lp)$: $N=L$. Hence $L$ is a two-sided ideal in $\U(\lp)$.
\end{proposition}
\begin{proof} The proof follows \cite{dixmier}, Section 7. There is a
  short exact sequence of rings
\[ 0\rightarrow N \rightarrow \U(\lp)\rightarrow K \rightarrow 0 .\]
Let $x_1,x_2,..,x_n$ be a basis with $\rho(x_1)\neq 0$ and
$\rho(x_2)=\cdots =\rho(x_n)=0$. Assume $x=x_1^{v_1}x_2^{v_2}\cdots
x_n^{v_n}\in \U(\lp)$. It follows $\rho_{\U}(x)=0$ if and only if
$v_2+\cdots +v_n\geq 1$. We get by the PBW Theorem a direct sum
decomposition
\[ \U(\lp)=\{x_1^{v_1}x_2^{v_2}\cdots x_n^{v_n}: v_2+\cdots +v_n\geq
1\}\oplus \{x_1^{v_1}:v_1 \geq 0\}.\]
It follows there is an inclusion of vector spaces
\[ \{x_1^{v_1}x_2^{v_2}\cdots x_n^{v_n}: v_2+\cdots +v_n\geq 1\} \subseteq N .
\]
One checks there is an equality of vector spaces
\[ x_1^{v_1}:v_1\geq 0\}=\{x_1^l\rho(x_1)^l: l\geq 0\} .\]
There is an inclusion
\[ \{x_1^{v_1}x_2^{v_2}\cdots x_n^{v_n}: v_2+\cdots +v_n\geq
1\}\oplus \{x_1^{v_1}:v_1 \geq 1 \}\subseteq L \]
hence $codim(L, \U(\lp))\leq 1$. Similarly $codim(N,\U(\lp))=1$. Since
$L\subseteq N$ it follows there is an equality $L=N$ and the Proposition is proved.
\end{proof}

\begin{definition} Let $char(\rho)=\U(\lg)\{x-\rho(x)\oneg: x\in
  \lp\}$ be the left \emph{character ideal} of $\rho$ in
  $\U(\lg)$. Let $char_l(\rho)=char(\rho)\cap \U_l(\lg)$ be the
  canonical filtration of $char(\rho)$.
\end{definition}

\begin{proposition} \label{char} The natural map $\phi:\U(\lg)\rightarrow M(\rho)$
  defined by $\phi(x)=x\otimes w$ defines an exact sequence of left
  $\U(\lg)$-modules
\[0\rightarrow char(\rho)\rightarrow \U(\lg)\rightarrow^\phi
M(\rho)\rightarrow 0.\]
\end{proposition}
\begin{proof} Let $X_{-\alpha_1},..,X_{-\alpha_m}$ be the basis for  $\ln$ constructed above and let
\[ X^P=X_{-\alpha_1}^{p_1}\cdots X_{-\alpha_m}^{p_m}\]
with $p_i\geq 0 $ integers. There is an isomorphism of right
$\U(\lp)$-modules
\[ \U(\lg)\cong \{X^P:p_i\geq 0\}\U(\lp) \]
hence we get an isomorphism of vector spaces
\[ M(\rho)\cong \{X^P\otimes w: p_i \geq 0\}.\]
Assume $X\in \U(\lg)$ is an element. We may write $X=\sum_P X^P x_P$
with $x_P\in \U(\lp)$. Assume
$\phi(X)=X\otimes w=0$. We get
\[ \phi(X)=X\otimes w=\sum_P X^Px_P\otimes w=\sum_P X^P\otimes
x_Pw=0\]
It follows
\[ x_Pw=0 \]
for all $P$ hence $x_P\in ker(\rho_{\U})=\U(\lp)\{y-\rho(y)\onep:y\in
\lp\}$. It follows $X=\sum_P X^Px_P\in char(\rho)$ hence
$ker(\phi)\subseteq char(\rho)$. One checks $char(\rho)\subseteq
ker(\phi)$ and the Proposition is proved.
\end{proof}

By Lemma \ref{basis} it follows 
\[ M(\rho)=K \{X^P\otimes w :p_i\geq 0\}.\]
Let 
\[ \ap=-p_1\alpha_1-\cdots -p_m\alpha_m\in \lh^*.\]
Assume $x\in \lh$.
We get
\[ x(X^P\otimes w)=[x,X^P]\otimes w+X^Px\otimes w=\]
\[\ap(x)X^P\otimes w+X^P\otimes xw=\]
\[\ap(x)X^P\otimes w+\lambda(x)X^P\otimes w=\]
\[(\lambda+\ap)(x)X^P\otimes w.\]
Hence
\[M(\rho)_{\lambda+\ap}=K\{X^P\otimes w\}.\]
It follows 
\[ M(\rho)=\oplus_P M(\rho)_{\lambda +\ap}.\]
Define
\[ M(\rho)_+=\oplus_{p_i \geq 1}M(\rho)_{\lambda+\ap}.\]
It follows 
\[ M(\rho)=M(\rho)_{\lambda}\oplus M(\rho)_+\]
where
\[ M(\rho)_{\lambda}=K\{\oneg \otimes w\}.\]
Let the vector $\oneg \otimes v\in M(\rho)$ be the \emph{canonical
  generator} of $M(\rho)$.

\begin{theorem} Let $M(\rho)$ be the scalar generalized Verma module
  associated to the character $\rho$. The following holds:
$M(\rho)$ contains a maximal non-trivial sub-$\lg$-module $K$.
The quotient $L(\rho)=M(\rho)/K$ is simple and $dim_K(L(\rho))\geq 2$.
Let $v=\overline{1\otimes w}\in L(\rho)$. The ideal $ann(v)$ is the
largest non-trivial left ideal in $\U(\lg)$ containing $char(\rho)$.
The vector $v$ satisfies the following: $\U(\lnp)v=0$. For all $x\in
\lh$ it follows $xv=\lambda(x)v$ hence $v$ has weight $\lambda$
\end{theorem}
\begin{proof} Consider the element $X_{-\alpha_m}\otimes w\in
  M(\rho)_+$ and look at the product
\[ X_{\alpha_m}X_{-\alpha_m}\otimes w.\]
We get
\[ X_{\alpha_m}X_{-\alpha_m}\otimes
w=[X_{\alpha_m},X_{-\alpha_m}]\otimes w+X_{-\alpha_m}\otimes
X_{\alpha_m}w=\]
\[ H_{\alpha_m}\otimes w=\oneg \otimes H_{\alpha_m}w=\]
\[\lambda(H_{\alpha_m})(\oneg \otimes w).\]
It follows 
\[ X_{\alpha_m}X_{-\alpha_m}\otimes w\in M(\rho)_{\lambda} \]
hence $M(\rho)_+$ is not $\lg$-stable. 

Assume $L\subsetneq M(\rho)$ is
a non-trivial $\lg$-stable module. It follows $L\cap
M(\rho)_{\lambda}=0$ hence $L\subseteq M(\rho)_+$. Let $K$ be the sum
of all non-trivial sub-$\lg$-modules of $M(\rho)$. It follows
$K\subseteq M(\rho)_+$ since $M(\rho)_+$ is not $\lg$-stable. Hence $K\subsetneq M(\rho)$ is a maximal
non-trivial sub-$\lg$-module of $M(\rho)$ and $L(\rho)=M(\rho)/K$ is
a simple quotient. It follows $dim_K(L(\rho))\geq 2$. 

One checks the
vector $v$ is annihilated by $\U(\lnp)$ and has weight $\lambda$. 

By Proposition \ref{char} there is an exact sequence of left
$\U(\lg)$-modules
\[ 0\rightarrow char(\rho)\rightarrow \U(\lg)\rightarrow
M(\rho)\rightarrow 0\]
hence there is an isomorphism
\[ M(\rho)\cong U(\lg)/char(\rho) \]
of left $\U(\lg)$-modules. It follows there is a bijection between the
set of left sub-$\lg$-modules of $M(\rho)$ and left sub-$\lg$-modules
of $\U(\lg)/char(\rho)$. This induce a bijection between the set of left
ideals in $\U(\lg)$ containing the ideal $char(\rho)$ and the set of
left sub-$\lg$-modules of $M(\rho)$. It follows the
submodule $K$ corresponds to a maximal non-trivial left ideal
$J$ in $\U(\lg)$ contaning $char(\rho)$. The ideal $J$ is by
definition the annihilator ideal of $v$: There is an equality
\[ J=ann(v).\]
The Theorem is proved. 
\end{proof}

\begin{corollary} The following holds:
\begin{align}
\label{C1}&\text{$L(\rho)$ is simple for all $\ul\in \mathbf{Z}^k$.}\\
\label{C2}&\text{If $\ul\neq \ul'$ it follows $L(\rho_{\ul})\neq L(\rho_{\ul'})$}\\
\label{C3}&\text{If $l_1,..,l_k<0$. It follows $dim_K(L(\rho))=\infty$}
\end{align}
\end{corollary}
\begin{proof}
Claim \ref{C1}: This is by definition of $L(\rho)$ since the submodule
$K$ is maximal.
Claim \ref{C2}: Assume $\phi:L(\rho_{\ul})\rightarrow L(\rho_{\ul'})$
is a map of $\lg$-modules. It follows $\phi$ is the zero map or an
isomorphism since the modules are simple. If it is an isomorphism it
follows the weights are equal. This implies $\ul=\ul'$ a
contradiction. The claim is proved.
We prove claim \ref{C3}:  Assume $dim_K(L(\rho))<\infty$. It follows $L(\rho)\cong
  \Vl$ where $\Vl$ has a highest weight vector $v'$ with highest
  weight $\lambda'=\sum_{j=1}^{n-1}l_j'\omega_j$ with $l_j'\geq
  0$. Since the  vector $v$ in $L(\rho)$ has weight $\lambda$ with
  $l_i<0$ one gets a contradiction. The Corollary follows.
\end{proof}

Assume $\lp=\lbp=\lh\oplus \lnp$. It follows from \cite{dixmier}, Section 7 the $\SL(E)$-module $L(\lambda+\delta)$ 
is isomorphic to $\Vl$ - the finite dimensional irreducible
$\SL(E)$-module with highest weight $\lambda=\sum_{i=1}^k
l_i\omega_{n_i}$. Here $l_i\geq 1$ is an integer for $i=1,..,k$. Hence
the class 
\[ \{ L(\rho_{\ul}): \rho_{\ul}:\lp\rightarrow K, \ul\in \mathbf{Z}^k \} \]
is a class of mutually non-isomorphic $\SL(E)$-modules parametrized by
a parabolic subalgebra $\lp\subseteq \sl(E)$ and a character $\rho:\lp
\rightarrow K$ containing the class of all finite dimensional irreducible $\SL(E)$-modules.

Let $K_l=\K\cap M_l(\rho)\subseteq M(\rho)$ be the induced filtration
on $K$. There is an exact sequence of $P$-modules
\[ 0\rightarrow K_l \rightarrow M_l(\rho)\rightarrow
L_l(\rho)\rightarrow 0  .\]

\begin{definition} Let $\{L_l(\rho)\}_{l\geq 0}$ be the
  \emph{canonical filtration} of $L(\rho)$. 
\end{definition}

It follows $\{L_l(\rho)\}_{l\geq 0}$ is a filtration of $L(\rho)$ by
$P$-modules.

\begin{corollary} Assume $L(\rho)=\Vl$ is a finite dimensional
  irreducible $\SL(E)$-module with highest weight vector $v$ and
  highest weight $\lambda$. It follows $L_l(\rho)\cong \U_l(\sl(E))v$.
\end{corollary}
\begin{proof} The proof is obvious.
\end{proof}

\section{Classical Verma modules and  annihilator ideals}

Let $\lbp=\lh\oplus \lnp$. It follows $\lbp$ is a sub Lie algebra of
$\sl(E)$. Let $(x,n), (y,m)$ be elements of $\lbp$. The Lie product on
$\sl(E)$ induce the following product on $\lbp$: define the following
action of $\lh$ on $\lnp$.
\[ ad:\lh\rightarrow \End(\lnp) \]
\[ ad(x)(n)=[x,n].\]
It follows $\lnp$ is a $\lh$-module. Let $ad(x)(h)=x(h)$.

\begin{lemma} Define
\[ [(x,n),(y,m)]=(0, x(m)-y(n)+[n,m]) \]
where $[,]$ is the bracket on $\lnp$. It follows the natural injection
$\lbp \rightarrow \sl(E)$ is a map of Lie algebras.
\end{lemma}
\begin{proof} The proof is obvious.
\end{proof}

Let $\mu\in \lh^*$ be a linear form on $\lh$. Let $L_\mu=Kw$ be the
free rank one $K$-module on $w$. Define the following map
\[ \tau_\mu:\lbp \rightarrow \End(L_\mu) \]
by
\[ \tau_\mu(x,n)(v)=\mu(x)v.\]

\begin{lemma} The map $\tau_\mu$ makes $L_\mu$ into a $\lbp$-module.
\end{lemma}
\begin{proof} It is clear
  $[\tau_\mu(x,n),\tau_\mu(y,m)]=\tau_\mu([(x,n),(y,m)])$
hence the claim is proved.
\end{proof}
It follows $L_\mu$ is a left $\lbp$-module. By definition
$\U(\sl(E))$ is a right $\lbp$-module and we may form the tensor
product
\[ M(\mu)=\U(\sl(E))\otimes_{\U(\lbp)}L_\mu.\]
It follows $M(\mu)$ is a left $G$-module.

\begin{definition} The $G$-module $M(\mu)$ is the \emph{Verma module}
  associated to the linear form $\mu \in \lh^*$.
\end{definition}

Let $\delta=\frac{1}{2}\sum_{\alpha\in R_+} \alpha$ and let $\lambda$
be an element in $\lh$. Let $L_{\lambda-\delta}=Kv$ be the free rank
one $K$-vectorspace on the element $v$.
By the result above we get a character
\[ \tau_{\lambda -\delta}:\lbp \rightarrow \End(L_{\lambda-\delta}) \]
defined by
\[ \tau_{\lambda-\delta}(h,n)=(\lambda-\delta)(h)\]
where $(h,n)$ is an element in $\lbp=\lh\oplus \lnp$. Let
$\alpha_1,..,\alpha_m$ be the $m$ distinct elements of $R_+$. Let
$X_{-\alpha_i}$ be an element in $\lg^{-\alpha_i}$. It follows the set
\[ X_{-\alpha_1},...,X_{-\alpha_m} \]
is a basis for $\lnm$ as vector space.
There is a decomposition $\lg=\lnm\oplus \lbp$ and it follows
$\U(\lg)$
is a free right $\U(\lbp)$-module as follows:
\[\U(\lg)=\{X^{p_1}_{-\alpha_1}\cdots X^{p_m}_{-\alpha_m}:
p_1,..,p_m\geq 0\}\U(\lbp).\]
\begin{lemma} Assume $X^P=X^{p_1}_{-\alpha_1}\cdots
  X^{p_m}_{-\alpha_m}$ with $p_1,..,p_m\geq 0$ integers. Let $x\in
  \lh$. It follows
\[ x(X^P)=(-p_1\alpha_1-\cdots -p_m\alpha_m)(x)X^P=\alpha_P(x)X^P .\]
\end{lemma}
\begin{proof} The proof is an easy calculation.
\end{proof}

Let $P=(p_1,..,p_m)$ with $p_i$ integers and let
$\alpha_P=p_1\alpha_1+\cdots +p_m\alpha_m$. Let
$X^P=X^{p_1}_{-\alpha_1}\cdots X^{p_m}_{-\alpha_m}$.
\begin{proposition} Assume $x\in \lh$. The following holds:
\begin{align}
\label{V1}&M(\lambda)=K\{X^P \otimes v: p_1,..,p_m\geq 0\}\\
\label{V2}&x(X^P \otimes v)=(\lambda-\delta-\alpha_P)(x)X^P \otimes v
\\
\label{V3}&M(\lambda)=\oplus_{p_1,..,p_m\geq 0}
M(\lambda)_{\lambda-\delta-\alpha_P}\\
\label{V4}&M(\lambda)_{\lambda-\delta}=K(1\otimes v)\\
\label{V5}&M(\lambda)_\mu=K\{X^P \otimes v:
\lambda-\delta-\alpha_P=\mu\}\\
\label{V6}&\U(\lnm)(1\otimes v)=M(\lambda)
\end{align}
\end{proposition}
\begin{proof}
We prove Claim \ref{V1}: There is a direct sum decomposition
\[ \lg=\lnm \oplus \lbp \]
and the set
\[ X_{-\alpha_1},..,X_{-\alpha_m} \]
is a basis for $\lnm$ as $K$-vector space. It follows by the
PBW-theorem that the set
\[ \{X^P : p_1,..,p_m\geq 0\} \]
is a basis for $\U(\lnm)$ as $K$-vector space. It follows $\U(\lg)$ is
isomorphic to
\[ \{X^P: p_1,..,p_m\geq \}\U(\lbp) \]
as free left $\U(\lbp)$-module. We get
\[ M(\lambda)=\U(\lg)\otimes_{\U(\lbp)}L_{\lambda-\delta}\cong \]
\[ \{X^P: p_1,..,p_m\geq
0\}\U(\lbp)\otimes_{\U(\lbp)}L_{\lambda-\delta} \cong \]
\[ K\{ X^P\otimes v: p_1,..,p_m\geq 0\} \]
and Claim \ref{V1} is proved.

We prove Claim \ref{V2}: By the previous Lemma it follows
for all $x\in \lh$
\[ [x,X^P]=\alpha_P(x)X^P.\]
We get
\[ x(X^P\otimes v)=[x,X^P]\otimes v+X^Px\otimes v=\]
\[ \alpha_P(x)X^P\otimes v+X^P\otimes (\lambda-\delta)(x)v=\]
\[ (\lambda-\delta-\alpha_P)(x)X^P \otimes v \]
and Claim \ref{V2} is proved.
We prove Claim \ref{V3}: By Claim \ref{V1} it follows $M(\lambda)$ has
a basis given as follows:
\[ \{ X^P\otimes v: p_1,..,p_m\geq 0\}.\]
Since for all $x\in \lh$ it follows
\[ x(X^P\otimes v)=(\lambda-\delta-\alpha_P)(x)X^P\otimes v \]
Claim \ref{V3} follows.

We prove Claim \ref{V4}: Let $X^P\otimes v\in M(\lambda)$
It follows $x(X^P\otimes v)=(\lambda-\delta-\alpha_P)(x)X^P \otimes
v$. Hence $X^P\otimes v$ is in $M(\lambda)_{\lambda-\delta}$ if and
only if
\[ \lambda -\delta -\alpha_P=\lambda-\delta.\]
Since $p_1,..,p_m\geq0$ this can only occur when $p_1=\cdots p_m=0$
hence $\alpha_P=0$. The Claim is proved.

Claim \ref{V5} is by definition.

We prove Claim \ref{V6}: Since the map
\[\phi: \U(\lnm)\rightarrow M(\lambda) \]
defined by
\[\phi(X^P)=X^P\otimes v\]
is an isomorphism of $\lnm$-modules Claim \ref{V6} follows. The
Proposition is proved.
\end{proof}
\begin{definition} The element $1\otimes v\in M(\lambda)$ is called
  the \emph{canonical generator} of $M(\lambda)$.
\end{definition}

\begin{proposition} There is an isomorphism of $\lnm$-modules
\[ \phi:\U(\lnm)\cong M(\lambda) \]
defined by
\[ \phi(X^P)=X^P \otimes v.\]
\end{proposition}
\begin{proof} Let $\lnm$ have basis $X_{-\alpha_1},..,X_{-\alpha_m}$
  and let $X^P=X_{-\alpha_1}^{p_1}\cdots X_{-\alpha_m}^{p_m}$. It
  follows $M(\lambda)$ has a basis consisting of elements $X^P\otimes
  v$. Hence the map $\phi$ is an isomorphism of vector spaces. One
  checks it is $\lnm$-linear and the Proposition follows.
\end{proof}

Let in the following $E=K\{e_1,..,e_n\}$ be an $n$-dimensional vector
space over $K$ and let $\lg=\sl(E)$ with $\lh$ in $\lg$ the abelian
sub-algebra of diagonal matrices.
Let $E_i=K\{e_1,..,e_i\}$ for $i=1,..,n-1$. It follows we get a
complete flag
\[ E_\bullet: 0\neq E_1\subseteq E_2\subseteq \cdots \subseteq
E_{n-1}\subseteq E
\]
in $E$ with $dim(E_i)=i$. Let $\underline{l}=(l_1,..,l_{n-1})$ with
$l_i\geq 0$ integers.
Let
\[W(\underline{l})=\sym^{l_1}(E)\otimes \sym^{l_2}(\wedge^2 E)\otimes
\cdots \otimes \sym^{l_{n-1}}(\wedge^{n-1} E) .\]
Let $v_i=\wedge^{i}E_i$ and let $v=v_1 {l_1}\otimes \cdots \otimes
v_{n-1}^{l_{n-1}}$ be a line in $W(\underline{l})$. Let $P$ in $SL(E)$
be the subgroup of elements fixing the flag $E_\bullet$. It follows
$P$ consists of upper triangular matrices with determinant one. Let
$\lp=Lie(P)$. It follows $\lp$ consists of upper triangular matrices
with trace zero.
Let
\[ \omega_i=L_1+\cdots +L_i \]
where
\[ \lh^*=K\{L_1,..,L_n\}/L_1+\cdots +l_n.\]
It follows $\omega_1,..,\omega_{n-1}$ are the fundamental weights for
$\lg$.
Let $\lambda=\sum_{i=1}^{n-1}l_i\omega_i$.
\begin{proposition}Let $x\in \lp$ be an element. The following formula
  holds:
\[x(v)=\sum_{i=1}^{n-1}l_i(a_{11}+\cdots +a_{ii})v \]
where $a_{ii}$ is the $i$'th diagonal element of $x$.
It follows the vector $v$ has weight $\lambda$.
The $\SL(E)$-module $\Vl$ generated by $v$ is an irreducible finite
dimensional $\SL(E)$-module with
highest weight vector $v$ and highest weight $\lambda$.
\end{proposition}
\begin{proof} The Proposition follows from \cite{dixmier}, Section 7
and an explicit calculation.
\end{proof}

Let $L_v$ in $\Vl$ be the line spanned by the vector $v$. It follows
the subgroup of $\SL(E)$ of elements fixing $L_v$ equals the group
$P$. We get a character
\[ \rho:\lp\rightarrow \End(L_v) \]
defined by
\[ \rho(x)v=x(v)=\sum_{i=1}^{n-1}l_i(a_{11}+\cdots +a_{ii})v.\]
We get an exact sequence of Lie algebras
\[ 0\rightarrow \lp_v \rightarrow \lp \rightarrow \End(L_v)
\rightarrow 0\]
where $\lp_v=ker(\rho)$.
\begin{definition} Let
\[ char(\rho)=\U(\lg)\{x-\rho(x)\oneg : x\in \lp\} \]
be the left \emph{character ideal} of $\rho$.
Let
\[ ann(v)=\{x\in \U(\lg): x(v)=0\} \]
be the left \emph{annihilator ideal} of $v$.
Let
\[ char_l(\rho)=char(\rho)\cap\U_l(\lg) \]
and
\[ ann_l(v)=ann(v)\cap\U_l(\lg)\]
for all $l\geq 1$.
\end{definition}

If $x\in \lh$ it follows $\rho(x)=\lambda(x)$ where $\lambda$ is the
highest weight of $\Vl$. By \cite{dixmier}, Section 7 the following
holds: Let for $\alpha \in B$ the integer $m_\alpha$ be defined as
follows:
\[ m_{\alpha}=\lambda(H_\alpha)+1 \]
where
\[ 0\neq H_\alpha \in [X_\alpha,X_{-\alpha}] .\]
If $\alpha_i=L_i-L_{i+1}$ it follows $X_{\alpha_i}=E_{i,i+1}$. It also
follows $X_{-\alpha_i}=E_{i+1,i}$. We get
\[ H_{\alpha_i}=E_{ii}-E_{i+1,i+1}.\]
\begin{lemma} The following holds for all $i=1,..,n-1$:
\[ m_{\alpha_i}=l_i+1.\]
\end{lemma}
\begin{proof} The proof is an easy calculation.
\end{proof}

We get by \cite{dixmier}, Proposisition 7.2.7, Section 7 the following description of the ideal $ann(v)$ in
$\U(\lg)$. Note that $\U(\lg)$ is a noetherian associative algebra and
the ideal $ann(v)$ is a left ideal in $\U(\lg)$. It follows $ann(v)$
has a finite set of generators. The following holds:
Let
\[ I(v)=\U(\lg)\lnp+\sum_{x\in   \lh}\U(\lg)(x-\lambda(x)\oneg).\]
It follows $I(v)\subseteq \U(\lg)$ is a left ideal.
It follows
\[ann(v)= I(v) +\sum_{\alpha\in B}\U(\lg)X^{m_\alpha}_{-\alpha}.\]
Let $I_l(v)=I(v)\cap \U_l(\lg)$ for all $l\geq 1$ integers.
Let for any element $x\in \U(\lg)$ 
\[fil(x)=min\{l : x\in \U_l(\lg) \} \]
be the \emph{filtration} of the element $x$. 

\begin{lemma} \label{capp} There is for every integer $l\ge 1$ an equality 
\[ ann_l(v)=I_l(v)+ \sum_{\alpha \in
  B}\U_{l-m_\alpha}(\lg)X^{m_\alpha}_{-\alpha} \]
of vector spaces.
\end{lemma}
\begin{proof} The inclusion 
\[ I_l(v)+ \sum_{\alpha \in
  B}\U_{l-m_\alpha}(\lg)X^{m_\alpha}_{-\alpha} \subseteq ann_l(v) \]
is obvious.
Assume 
\[ x=u+v\in ann_l(v)=(I(v)+\sum_{\alpha\in
  B}\U(\lg)X^{m_\alpha}_{-\alpha})\cap \U_l(\lg) \]
with
\[ u \in I(v) \]
and
\[ v\in \sum_{\alpha\in B}\U(\lg)X^{m_\alpha}_{-\alpha}.\]
It follows $fil(u),fil(v)\leq fil(x)$. Hence $u,v\in \U_l(\lg)$. 
It follows $u\in I_l(v)$ and 
\[ v\in \sum_{\alpha\in B}U_{l-m_\alpha}(\lg)X^{m_\alpha}_{-\alpha}.\]
The Lemma follows.
\end{proof}

\begin{proposition} \label{characterideal} There is an equality
\[ I(v)=char(\rho) \]
of left ideals in $\U(\lg)$.
\end{proposition}
\begin{proof} Recall there is a triangular decomposition
  $\lg=\lnm\oplus \lh \oplus \lnp$.
It follows there is an equality  $\lp=\lh\oplus \lnp=\lbp$ of Lie
algebras.
Assume $x\in \lnp\subseteq \lp$. It follows $\rho(x)=0$, hence
$x(v)=\rho(x)v=0$. Assume $x\in \lh\subseteq \lp$. It follows
$\rho(x)=\lambda(x)$ where $\lambda $ is the highest weigt of $\Vl$
with
\[ \lambda =\sum_{i=1}^{n-1}l_i\omega_i.\]
It follows there is an element $x\in \lh$ with $\lambda(x)\neq
0$. Moreover we may choose elements $h_1,..,h_{n-2}$ in $\lh$ with
$\rho(h_i)=\lambda(h_i)=0$ and
\[ \{x,h_1,..,h_{n-2}\} \]
is a basis for $\lh$. Assume $y(x-\rho(x)\oneg)\in char(\rho)$. It
follows
\[ x\in \lp=\lh\oplus \lnp.\]
It follows $x=x_1+x_2$ with $x_1\in \lh$ and $x_2\in \lnp$.
It follows
\[y(x-\rho(x)\oneg)=y(x_1+x_2-\rho(x_1+x_2)\oneg)=y(x_1-\rho(x_1)\oneg)+y(x_2-\rho(x_2)\oneg).\]
Since $x_1\in \lh$ it follows
\[y(x_1-\rho(x_1)\oneg)=y(x_1-\lambda(x_1)\oneg) \in \sum_{x\in
  \lh}\U(\lg)(x-\lambda(x)\oneg) \subseteq I(v).\]
Since $x_2\in \lnp$ it follows $\rho(x_2)=0$. We get
\[ y(x_2-\rho(x_2)\oneg)=yx_2\in \U(\lg)\lnp.\]
It follows $y(x-\rho(x)\oneg)\in I(v)$ and
we have proved the inclusion $char(\rho)\subseteq I(v)$.
We prove the reverse inclusion: Assume $x=x_1+x_2\in I(v)$ with
$x_1=yn\in \U(\lg)\lnp$ and $y\in \U(\lg)$ and $n\in \lnp$. Assume
also $x_2=z(u-\lambda(u)\oneg)$ with $z\in \U(\lg)$ and $u\in \lh$. It
follows $\rho(n)=0$ and $x_1=y(n-\rho(n)\oneg)\in char(\rho)$. Also
$x_2=z(u-\lambda(u)\oneg)=z(u-\rho(u)\oneg)\in char(\rho)$. It follows
$x=x_1+x_2\in char(\rho)$ and we have proved the reverse inclusion
$I(v)\subseteq char(\rho)$. The Proposition is proved.
\end{proof}

\begin{corollary} \label{characterl} There is for every $l\geq 1$ an
  equality
\[ I_l(v)=char_l(\rho) \]
of filtrations.
\end{corollary}
\begin{proof} The Corollary follows from Proposition
  \ref{characterideal}.
\end{proof}

Recall $\lambda=\sum_{i=1}^{n-1}l_i\omega_i$. Define
$m(\lambda)=min_{i=1,..,n-1}\{l_i\}$.

\begin{theorem} \label{equality} For all $1\leq l \leq m(\lambda)$
  there is an equality
\[ ann_l(v)=char_l(\rho) \]
of filtrations.
\end{theorem}
\begin{proof} By Lemma \ref{capp} there is an eqality
\[ ann_l(v)=I_l(v)+\sum_{\alpha_i \in
  B}\U_{l-m_{\alpha_i}}(\lg)X^{m_{\alpha_i}}_{-\alpha_i}\]
of vector spaces.
We get
\[ X^{m_{\alpha_i}}_{-\alpha_i}=E_{i+1,i}^{l_i+1}\]
for $i=1,..,n-1$. We get from Corollary \ref{characterl}
\[
ann_l(v)=I_l(v)+\sum_{i=1}^{n-1}\U_{l-l_1-1}(\lg)E_{i+1,i}^{l_i+1}=\]
\[ char_l(\rho)+ \sum_{i=1}^{n-1}\U_{l-l_1-1}(\lg)E_{i+1,i}^{l_i+1}.\]
If $1\leq l \leq m(\lambda)$ it follows $l\leq l_i$ for all
$i=1,..,n-1$.
It follows $l-l_i-1<0$ hence we get an equality
\[ ann_l(v)=char_l(\rho) \]
and the claim of the Theorem follows.
\end{proof}

Let $\lp\subseteq \sl(E)$ be a parabolic subalgebra and let $\ln$ be
the complementary Lie algebra as defined in \cite{maa1}. It follows
there is a direct sum decomposition as vector spaces $\lp\oplus
\ln\cong \sl(E)$. We get inclusions $\ln\subseteq \lnm$ and $\lbp
\subseteq \lp$ of Lie algebras. Let $\rho:\lp\rightarrow
\End(L_\rho)$ with induced character
\[ \tilde{\rho}:\lbp\rightarrow \End(L_{\tilde{\rho}}).\]

\begin{lemma} There is a surjective map of left $\U(\sl(E))$-modules
\[\phi:M(\tilde{\rho})\rightarrow M(\rho) .\]
It follows $ker(\phi)=char(\rho)/char(\tilde{\rho})$.
\end{lemma}
\begin{proof} There is a natural map
\[ f:\U(\sl(E))\times L_\rho \rightarrow
\U(\sl(E))\otimes_{\U(\lp)}L_\rho \]
defined by
\[ f(z,w)=z\otimes w.\]
This map induce a surjection
\[ \phi:U(\sl(E))\otimes_{\U(\lbp)}L_\rho \rightarrow
\U(\sl(E))\otimes_{\U(\lp)}L_\rho \]
of left $\U(\sl(E))$-modules and the first claim is proved. The second
claim follows easily and the Lemma is proved.
\end{proof}

Let $K\subseteq M(\rho)$ be the unique maximal $\sl(E)$-stable
submodule and let $K'\subseteq M(\tilde{\rho})$ be the unique
$\sl(E)$-stable submodule. It follows $\phi(K')\subseteq K$ since the
image of a $\sl(E)$-stable module is $\sl(E)$-stable. We get a
commutative diagram of maps of $\sl(E)$-modules where the middle
vertical arrow is surjective:
\[
\diagram 0 \rto  & K' \rto \dto & M(\tilde{\rho}) \rto^{p'} \dto^{\phi} & L(\tilde{\rho}) \dto^{\phi'} \rto & 0 \\
0 \rto  & K \rto  & M(\rho) \rto^{p} & L(\rho)  \rto & 0 \\
\enddiagram
.\]

\begin{lemma} The map $\phi'$ is an isomorphism of $\sl(E)$-modules.
\end{lemma}
\begin{proof} Since the map $\phi$ is a surjection of $\sl(E)$-modules
  it follows $\phi'$ is a surjection. Since the modules $L(\rho)$ and
  $L(\tilde{\rho})$ are simple it follows $ker(\phi')$ is zero. The
  Lemma is proved.
\end{proof}

Hence we get no new simple $\SL(E)$-modules when we consider simple quotients $L(\rho)$
where $\lp$ is a parabolic subalgebra of $\sl(E)$ with $\lp\neq
\lbp$. The modules $L(\rho)$ may be constructed using simple quotients of
classical Verma modules as done in \cite{dixmier}.

Let $u=1\otimes w\in M(\rho)$ be the canonical generator and let $v$
be its image in $L(\rho)$ under the canonical projection map
$M(\rho)\rightarrow L(\rho)$. Let $ann(v)$ be the left annihilator
ideal in $\U(\sl(E))$ of the vector $v$.

Recall the exact sequence
\[ 0\rightarrow K_l \rightarrow M_l(\rho) \rightarrow L_l(\rho)
\rightarrow 0\]
from the previous section.

\begin{proposition} \label{zero} The following holds:
\[ K_l=0 \text{ if and only if }ann_l(v)=char_l(\rho).\]
\end{proposition}
\begin{proof} There is a commutative diagram of exact sequences of
  left $\sl(E)$-modules
\[
\diagram 0 \rto & char(\rho) \rto \dto & \U(\sl(E)) \rto^\phi \dto &
M(\rho) \rto \dto & 0 \\
0 \rto & ann(v) \rto & \U(\sl(E)) \rto^{\tilde{\phi}} & L(\rho) \rto & 0
\enddiagram
\]
One checks that $\phi(ann(v))=K$ and $\phi(ann_l(v))=K_l$ for all
$l\geq 1$. Assume $ann_l(v)=char_l(\rho)$. It follows
\[ K_l=\phi(ann_l(v))=\phi(char_l(\rho))=0 \]
since the diagram above is exact. Assume $K_l=0$ and let $x\in
ann_l(v)$. It follows
\[ \phi(x)=x\otimes w\in K_l\]
hence $\phi(x)=0$. It follows $x\in char_l(\rho)$ since the diagram
has exact rows. It follows $ann_l(v)=char_l(\rho)$ and the Proposition follows.
\end{proof}

Assume $\lp=\lbp$ and $\rho|_\lh=\lambda$ with $\lambda=\sum_i l_i
\omega_i$ and $l_i\geq 0$ for all $i$.

\begin{corollary} \label{iso} If $1\leq l \leq m(\lambda)$ it follows there is an
  isomorphism
\[ M_l(\rho)\cong L_l(\rho) \]
of left $\SL(E)$-modules.
\end{corollary}
\begin{proof} The Corollary follows from Proposition \ref{zero} and Theorem \ref{equality}.
\end{proof}

Let $P\subseteq \SL(E)$ be the parabolic subgroup with $Lie(\lp)=P$
and let $\L\in \Pic^{\SL(E)}(\SL(E)/P)$ be the linebundle with
$\H^0(\SL(E)/P, \L)^*=L(\rho)=\Vl$. Let $X=\SL(E)/P$ and let
$\Pr^l_X(\L)$ be the $l$'th order jetbundle of $\L$ as defined in \cite{maa0}.

\begin{corollary} \label{jets} For all $1\leq l \leq m(\lambda)$ it follows there
  is an isomorphism
\[ \Pr^l_{X}(\L)(e)^*\cong \U_l(\sl(E))\otimes_{\U(\lp)}L_\rho \]
of $P$-modules.
\end{corollary}
\begin{proof} By \cite{maa0}, Theorem 3.10 there is an isomorphism
\[ \Pr^l_X(\L)(e)^*\cong \U_l(\sl(E))v\cong L_l(\rho) \]
of $P$-modules when $1\leq l \leq m(\lambda)$. The Corollary now
follows from Corollary \ref{iso}.
\end{proof}

Corollary \ref{jets} gives an algebraic proof of a result proved in
\cite{smoke} over the complex numbers.


\begin{thebibliography}{4}

\bibitem{dixmier} J. Dixmier, Enveloping Algebras, Graduate Studies in
  Mathematics,\emph{American Math. Society} (1996)

\bibitem{maa0} H. Maakestad, On jetbundles and generalized Verma
  modules II, \emph{Preprint math arXiv 0903.3291}

\bibitem{maa1} H. Maakestad, On the annihilator ideal of a highest
  weight vector, \emph{Preprint math arXiv:1003.3522}

\bibitem{smoke} W. Smoke, Invariant differential operators,
  \emph{Trans. Am. Math. Soc.} no. 127 (1967)

\end{thebibliography}
\end{document}